\RequirePackage{fix-cm}
\documentclass[smallextended]{svjour3}       
\smartqed  
\usepackage[nointegrals]{wasysym}
\usepackage[misc]{ifsym}
\usepackage{mathrsfs,amsmath,amssymb,amsfonts,enumerate,lipsum,appendix}
\usepackage{graphicx}
\usepackage[dvipsnames]{xcolor}
\definecolor{ao}{rgb}{0.0, 0.5, 0.0}
\definecolor{lasallegreen}{rgb}{0.03, 0.47, 0.19}
\usepackage{hyperref}
\hypersetup{
	colorlinks=true,
	linkcolor=red,
	filecolor=magenta,      
	urlcolor=cyan,
	citecolor=lasallegreen,
}
\usepackage[margin=0.75in]{geometry}

\let\oldnorm\norm
\def\norm{\@ifstar{\oldnorm}{\oldnorm*}}
\newcommand{\md}[1]{\left\vert #1 \right\vert}
\newcommand{\dq}{\mathrel{\mathop:}=}

\newcommand{\RN}{\mathbb{R}^N}
\newcommand{\pr}[1]{{#1}^{\prime}}
\newcommand{\al} {\alpha}
\newcommand{\pa} {\partial}

\newcommand{\de} {\delta}
\newcommand{\De} {\Delta}
\newcommand{\Dep} {\Delta_p}
\newcommand{\ga} {\gamma}
\newcommand{\Ga} {\Gamma}
\newcommand{\om} {\omega}
\newcommand{\Om} {\Omega}
\newcommand{\la} {\lambda}
\newcommand{\La} {\Lambda}

\newcommand{\Gr} {\nabla}

\newcommand{\noi} {\noindent}

\newcommand{\dx}{{\,\rm d}x}
\newcommand{\dt}{{\,\rm d}t}
\newcommand{\dr}{{\,\rm d}r}
\newcommand{\ds}{{\,\rm d}s}
\newcommand{\da}{{\,\rm d}\al}
\newcommand{\dd}{{\,\rm d}\de}

\newcommand{\dsig}{{\,\rm d}\sigma}

\newcommand{\A}{{\mathcal A}}
\newcommand{\C}{{\mathcal C}}
\newcommand{\E}{{\mathcal E}}
\newcommand{\F}{{\mathcal F}}
\newcommand{\G}{{\mathcal G}}
\newcommand{\M}{{\mathcal M}}
\newcommand{\hs}{\hspace}
\newcommand{\R}{{\mathbb R}}
\newcommand{\N}{{\mathbb N}}

\def\S{\mathcal{S}}
\def\wp{{{W}^{1,p}}(\Om)}

\usepackage[hyperpageref]{backref}
\begin{document}
\title{On reverse Faber-Krahn inequalities 
}

\titlerunning{Reverse Faber-Krahn inequalities}        

\author{T. V. Anoop, K. Ashok Kumar 
}

\authorrunning{T. V. Anoop \and K. Ashok Kumar} 

\institute{T.V. Anoop (\Letter)\at
              Department of Mathematics, Indian Institute of Technology, Chennai 600036, India \\
              \email{anoop@iitm.ac.in}           
           \and
           K. Ashok Kumar \at
              Department of Mathematics, Indian Institute of Technology, Chennai 600036, India\\
              \email{srasoku@gmail.com}
}

\date{Received: date / Accepted: date}

\maketitle

\begin{abstract} 
	Payne-Weinberger showed that `{\it among the class of membranes with given area $A$, free along the interior boundaries and fixed along the outer boundary of given length $L_0$, the  annulus $\Omega^\#$ has the highest fundamental frequency,}' where $\Omega^\#$ is a concentric annulus with the same area as $\Omega$ and the same outer boundary length as $L_0$. We extend this result for the higher dimensional domains and $p$-Laplacian with $p\in (1,\infty),$ under the additional assumption that the outer boundary is a sphere.
	As an application, we prove that the nodal set of the second eigenfunctions of $p$-Laplacian (with mixed boundary conditions) on a ball and a concentric annulus  cannot be a concentric sphere.
\end{abstract}
\section{Introduction}\label{section 1}
Let us first recall the famous conjecture by Lord Rayleigh from his book `The theory of sound' \cite{Rayleigh} published in 1877. He conjectured that   
{\it `among all  planar domains $\Om$ of fixed area, the disk is the  domain that minimises the first Dirichlet eigenvalue $\la_1(\Om)$ of the Laplacian.'} This conjecture was open for a very long time. The first proof of this conjecture was published in 1923 by Faber \cite{Faber}. In 1925, Krahn \cite{Krahn}  gave an independent proof for this conjecture, and later he extended the same to the higher dimension.  Now this result is collectively known as the Faber-Krahn inequality, and it states that: 
\begin{equation}\label{F-K}
\la_1(\Om^*)\leq \la_1(\Om)
\end{equation} 
{\it with the equality if and only if $\Om$  is a ball (up to a measure zero set), where $\Om^*$ is a ball of the same measure as $\Om$.}  For more on Faber-Krahn inequality and related results, we refer to \cite{Payne} and \cite{Polya-Szego}.

Among the class of domains having a fixed volume, $\la_1(\Om)$ is not bounded above. For example, consider a class of arbitrary long rectangular domains with the same volume (see Remark 3.13 of \cite{Anoop}). However, if we fix both the surface area and the volume one can give upper bounds for $\la_1(\Om )$, see \cite{Makai1} and \cite{Polya}. In \cite{P-W} Payne-Weinberger considered  this problem for bounded planar membrane with fixed  outer boundary having the interior boundaries along which it is free. They showed that `{\it among the class of membranes with given area $A$ and given length $L_0$ of the outer boundary,  the annulus $\Om^\#$ has the highest fundamental frequency,}' where $\Om^\#$ is a concentric annulus with the same area as $\Om$ and the same outer boundary length as $L_0$. In other words, this result can be stated in the form of an inequality as below:
\begin{equation}\label{Rev1}
\nu_1(\Om) \le  \nu_1(\Om^\#),
\end{equation} where $\nu_1(\Om),\, \nu_1(\Om^\#)$ are first eigenvalues of the following mixed eigenvalue problems:\\
\begin{minipage}[t]{0.49\textwidth}
	\begin{equation*}
	\left.
	\hs{-0.8cm}
	\begin{aligned}
	-\De u&= \nu u \quad \text{in } \, \Om,\\
	u&=0 \text{ on }  \Ga_0,\\
	\frac{\pa{u}}{\pa \eta}&=0  \text{ on }   \pa\Om\setminus\Ga_0;
	\end{aligned}
	\right\}
	\end{equation*}\hspace*{3cm}
\end{minipage}
\begin{minipage}[t]{0.49\textwidth}
	\begin{equation*}
	\left.
	\hs{-0.8cm}
	\begin{aligned}
	-\De u&= \nu u \quad \text{in } \, \Om^\#,\\
	u&=0 \text{ on }  \Ga_0^\#,\\
	\frac{\pa{u}}{\pa \eta}&=0  \text{ on }   \pa\Om^\#\setminus \Ga_0^\#,
	\end{aligned}
	\right\}
	\end{equation*}
\end{minipage}\\
where $\Ga_0$ and  $\Ga_0^\#$ are the outer boundaries of $\Om$ and $\Om^\#$ respectively.

We call the inequality in \eqref{Rev1}, as the reverse Faber-Krahn (R-F-K) inequality. A  similar inequality for the second Neumann eigenvalue (the first non-zero eigenvalue), namely $\mu_2(\Om)\le \mu_2(\Om^*)$ is obtained in  \cite{Szego} by Szeg\"o  for planar domains, and  in \cite{Weinberger-Neumann} for  higher dimensions by Weinberger. The proof of Szeg\"o mainly rely on the conformal mapping, and the proof of Weinberger based on construction of a test function using a radial function together with a suitable translation of the origin. In \cite{P-W}, for  proving the inequality \eqref{Rev1}, authors used the interior parallels and an isoperimetric inequality, deduced from an inequality for the  interior parallels due to  B. Sz. Nagy \cite{Nagy}. Next, we briefly describe the interior parallels and Nagy's inequality.

\noi{\bf The parallel sets:}  Let $\Om$ be a bounded domain in $\R^2$ and $\Ga_0$ be its outer boundary. Let $L(\de)$   be the measure of the inner parallel to $\Ga_0$-{\it  the set of all points that are inside $\Ga_0$ and at a distance $\de$ from $\Ga_0$.} Nagy showed that the function $L(\de)$ is defined for almost every $\de$ and  it satisfies 
\begin{equation}\label{Ineq:Nagy}
L(\de) \le |\Ga_0|-2\pi \de.
\end{equation}
Similarly, if $L(\de)$ is the measure of the outer parallel curve to $\Ga_1$, an interior boundary; he showed that $L(\de)$ satisfies 
\begin{equation}\label{Ineq2:Nagy}
L(\de) \le |\Ga_1|+2\pi \de.
\end{equation}

\noi Recall the Steiner formula $L(\de)=|\Ga_1|+2\pi \de$ for the convex planar domains, (see Chapter 4 of \cite{Schneider}). Indeed, Nagy's second inequality is an extension of Steiner formula for the non-convex planar domains. 

The {\it interior parallels} are the part of the parallels that lies in $\Om.$ If $l(\de)$ is the measure of the interior parallel at a distance $\de$, then $l(\de)\le L(\de).$ The interior parallels and Nagy's inequalities were used in \cite{Makai2}, \cite{Makai1} and \cite{Polya} for finding the upper bounds for $\la_1(\Om)$ among the class of domains having a fixed volume and surface area. In \cite{P-W}, using Nagy's inequality \eqref{Ineq:Nagy},  authors derived the following isoperimetric inequality:
\begin{equation}\label{P-W ineq}
l(\de)^2\leq |\Ga_0|^2-4\pi v(\de),
\end{equation} 
\noi where $v(\de)=\int_0^\de l(s)d\ds.$ Further, using  $v(\de)$  and the first eigenfunction of the concentric domain $\Om^\#$, they have constructed a test function on $\Om$ whose level sets coincide with interior parallels and the Rayleigh quotient is smaller than $\nu_1(\Om^\#).$  In \cite{Hersch}, Hersch gave another proof for \eqref{Rev1} using a class of test functions known as `web functions' (the test function that depends only on the distance from $\Ga_0$),  for more on web functions see \cite{BCT} and \cite{Gazzola} and the references therein. 

Hersch also considered the mixed eigenvalue problem with the Dirichlet condition  on an interior boundary and with the Neumann condition on the rest of the boundaries. By an ingenious parametrization $t(\de)=\int_0^\de \frac{\ds}{l(s)}$, he related this problem with the fundamental frequency of a vibrating string of infinite length (with varying mass distribution) fixed at one end. Then he used Nagy's inequality \eqref{Ineq2:Nagy} to show that `{\it among the class of membranes with given area $A$ and given length $L_1$ of an inner boundary $\Ga_1$, the annulus $\Om_\#$ has the maximal fundamental frequency,'} where $\Om_\#$ is a concentric annulus with the same area as $\Om$ and the same inner boundary length as $L_1$. As before, this can be restated as a R-F-K inequality:
\begin{equation}\label{Rev2}
\tau_1 (\Om) \le  \tau_1(\Om_\#),
\end{equation} where $\tau_1(\Om),\, \tau_1(\Om_\#)$ are first eigenvalues of the following problems:\\
\begin{minipage}[t]{0.49\textwidth}
	\begin{equation*}
	\left.
	\hs{-0.8cm}
	\begin{aligned}
	-\De u&= \tau u \quad \text{in } \, \Om,\\
	u&=0 \text{ on }  \Ga_1,\\
	\frac{\pa{u}}{\pa \eta}&=0  \text{ on }   \pa\Om\setminus\Ga_1;
	\end{aligned}
	\right\}
	\end{equation*}\hspace*{3cm}
\end{minipage}
\begin{minipage}[t]{0.49\textwidth}
	\begin{equation*}
	\left.
	\hs{-0.8cm}
	\begin{aligned}
	-\De u&= \tau u \quad \text{in } \, \Om_\#,\\
	u&=0 \text{ on }  \Ga_1^\#,\\
	\frac{\pa{u}}{\pa \eta}&=0  \text{ on }   \pa\Om_\#\setminus \Ga_1^\#.
	\end{aligned}
	\right\}
	\end{equation*}
\end{minipage}

\noi The mixed eigenvalue problem with the Dirichlet condition on both $\Ga_1,\Ga_0$, and with Neumann condition on the rest of the boundaries also studied in \cite{Hersch}. In this case, by an effective use of  the ``effectless cut'' due to Weinberger \cite{Weinberger} together with \eqref{Rev1} and \eqref{Rev2}, he showed that `{\it among the class of membranes with given area $A$, given length $L_0$ of the inner boundary $\Ga_0$ and given length $L_1$ of the outer boundary $\Ga_1$ satisfying $L_1^2- L_0^2=4\pi A$, the annulus has the highest fundamental frequency $\la_1$}.' In particular, an annulus $\Om=B_1\setminus \overline{B_0} \subset \R^2$  satisfies this relation. Thus by  Hersch's result:
\begin{equation}\label{Rev3}
\la_1 (B_1\setminus \overline{B_0}) \le  \la_1(B_1^*\setminus \overline{B_0^*}).
\end{equation}

In \cite{Ramm-Shivakumar}, Ramm and Shivakumar  conjectured the inequality  \eqref{Rev3} for $N\ge 3$,  with a numerical justification. In fact, they have conjectured a stronger result - \textit {`$\la_1(B_1\setminus \overline{B_0})$ strictly decreases when the inner ball moves towards the outer boundary.'} An  analytic  proof  for this  conjecture using an argument of M. Ashbaugh was   published later in an arxiv paper (arxiv:math-ph/9911040) by the same authors.  At the same time,    Harrell et al. \cite{HarrelKroger} and Kesavan  \cite{Kesavan} independently proved the strict monotonicity of $\la_1(B_1\setminus\overline{B_0})$. All the  proofs are based mainly  on the Hadamard perturbation formula. In  \cite{Anisa2}, the authors  studied   the monotonicity   of   the first Dirichlet eigenvalue of $p$-Laplacian, defined by $\Dep u =\mbox{div}(|\nabla u|^{p-2}\nabla u)$ for $p\in(1,\infty)$. They showed that \textit {`$\la_1(B_1\setminus \overline{B}_0)$  decreases when the inner ball $B_0$ moves towards the boundary of the outer ball $B_1$.'} Indeed, this result imply \eqref{Rev3} for the $p$-Laplacian. The strict monotonicity  is obtained in \cite{Anoop} and this result asserts that the equality happens in \eqref{Rev3} only if the balls $B_1$ and $B_0$ are concentric.

In this article, we study the R-F-K inequalities  \eqref{Rev1} and \eqref{Rev2} for the $p$-Laplacian with $p\in(1,\infty)$ and $N\ge 2$. We consider $\mathcal{C}^1$-smooth multiply connected bounded domains in $\R^N$.  For such a domain $\Om$,  we denote the outer boundary by $\Ga_0$ and the boundaries of the interior holes by $\Ga_1, \ga_1,\ldots \ga_n$. We use the symbol $|\,\cdot\,|$ for both $N$ dimensional volume measure and the $(N-1)$-dimensional surface measure. The Lebesgue measure of the unit ball in $\RN$ is denoted by $\om_N.$
Further, we set 

\begin{enumerate}[align=left,]
	\item[$\Om^*=$]  the  ball centred at the origin  with $|\Om^*|=|\Om|,$
	\item[$\Om^\#=$] the concentric annulus centred at the origin  with the same volume as $\Om$ and the outer surface measure as $\Ga_0,$
	\item[$\Om_\#=$] the concentric annulus centred at the origin  with the same volume as $\Om$
	and the inner surface measure as$\Ga_1.$
\end{enumerate}

In general, for a fixed domain $\Om$ the sets $\Om^\#$ and $\Om_\#$ are different. However, they coincide if $\Om$ satisfies the relation $|\Ga_0|^{\pr N}-|\Ga_1|^{\pr N}=C(N)|\Om|,$ where  $N^{\prime}$ is the H\"{o}lder conjugate of $N$ and $C(N)= N^{N^{\prime}} \om_N^{N^{\prime}-1}.$ We mainly consider two types of  mixed eigenvalue problems:  
 \begin{enumerate}[(i)]
\item the Dirichlet condition  on the outer boundary $\Ga_0$,
\item the Dirichlet condition  on an inner boundary $\Ga_1$,
\end{enumerate}
and the Neumann condition  on the rest of the boundaries. More precisely, for $p\in(1,\infty)$ we consider the following eigenvalue problems:

\begin{minipage}[t]{0.47\textwidth}
	\begin{equation}\label{ND}\tag{N-D}
	\left.
	\hs{-0.8cm}
		\begin{aligned}
	 -\De_p u&= \nu |u|^{p-2}u \quad \text{in } \, \Om,\\
	 u&=0 \text{ on }  \Ga_0,\\
	 \frac{\pa{u}}{\pa \eta}&=0  \text{ on }   \Ga_1\cup\left(\cup_{i=1}^n\ga_i\right);
	\end{aligned}
	\right\}
	\end{equation}\hspace*{3cm}
	\end{minipage}
	\begin{minipage}[t]{0.47\textwidth}
	\begin{equation}\label{DN}\tag{D-N}
	\left.
	\hs{-0.8cm}
		\begin{aligned}
	 -\De_p u&= \tau |u|^{p-2}u \quad \text{in } \, \Om,\\
	 u&=0 \text{ on }  \Ga_1,\\
	 \frac{\pa{u}}{\pa \eta}&=0  \text{ on }   \Ga_0\cup\left(\cup_{i=1}^n\ga_i\right).
	\end{aligned}
	\right\}
	\end{equation}
	\end{minipage}
	
	 \noi For $i=0,1,$ let 	$W_{\Ga_i}\dq\left\{u\in \wp : u|_{\Ga_i}=0 \right\}.$
	We say a real number $\nu$ is an eigenvalue of \eqref{ND}, if there exists $u\in W_{\Ga_0}\setminus\{0\}$ such that
	\begin{equation*}
	  \int_\Om |\nabla u|^{p-2}\nabla u\cdot \nabla v \dx = \nu \int_\Om |u|^{p-2}uv\dx,\, \forall v\in W_{\Ga_0}.  
	\end{equation*}
\noi	Similarly a real number $\tau$ is an eigenvalue of \eqref{DN}, if there exists $u\in W_{\Ga_1}\setminus\{0\}$ such that 
	\[\int_\Om |\nabla u|^{p-2}\nabla u\cdot \nabla v \dx = \tau \int_\Om |u|^{p-2}uv\dx,\, \forall v\in W_{\Ga_0}. \]
	The corresponding nonzero solutions are called the eigenfunctions associated to these eigenvalues. 
	Now consider
	\[J(u)\dq{\int_{\Om}\md{\nabla u}^p}\dx, \mbox{ and }\S\dq\left\{v\in \wp:\int_\Om |u|^p\dx=1\right\}.\]
	 For $i=0,1,$ one can easily see that the critical values of $J$  on $\S\cap W_{\Ga_i}$  are precisely the eigenvalues of \eqref{ND} and \eqref{DN}  respectively.  Thus  the classical Ljusternik-Schnirrelman theory, ensures
	 the existence of infinitely many eigenvalues  for both \eqref{ND} and \eqref{DN}, see Proposition \ref{propo:appendix 1}. 
	  In particular, the first eigenvalues  have the following variational characterisation:
	\[\nu_1(\Om)=\inf\limits_{u\in \S\cap W_{\Ga_0}} J(u),\ \tau_1(\Om)=\inf\limits_{u\in\S\cap W_{\Ga_1}} J(u).\]
	Both $\nu_1(\Om)$ and $\tau_1(\Om)$  are simple and the first eigenfunctions are of constant sign.  Moreover, the eigenfunctions corresponding to $\nu_1(\Om^\#)$  and $\tau_1(\Om_\#)$ are radial (see Appendix \ref{Appendix A}).
	
  The R-F-K inequalities, except \eqref{Rev3} for the annular regions, unfortunately did not get any attention for the domains in the higher dimensions and also for  the operators different from the Laplacian. Unlike in the case of Faber-Krahn inequality,  the uniqueness of the domains for which the equality holds in \eqref{Rev1} (similarly in \eqref{Rev2}) are not well understood, even for the planar domains. For $N=2,$ we extend the results of Payne-Weinberger and Hersch for the $p$-Laplacian. For $N\ge 3$, under an additional assumption on the boundary of $\Om,$ we  prove the R-F-K inequalities. To state this additional assumption, we make the following definition:

 \noi{\bf Definition:}
 Let $A,B$ be measurable sets in $\RN$. We say that $A$ is a $\mu$-translate of $B$, if there exists $x\in \R^N$ such that $|(A+x)\bigtriangleup B|=0,$ where $\bigtriangleup$ is the  symmetric difference of sets.

Now we state the R-F-K inequality for $\nu_1$: 
\begin{theorem}\label{thm 1}
Let $p\in(1,\infty)$ and $\Om\subset \R^N$ with $N\ge 3$ be a $\C^1$-smooth multiply connected bounded domain. If the outer boundary $\Ga_0$ of $\Om$ be  a sphere, then 
\begin{enumerate}[(i)]
 \item $ \nu_1(\Om) \le  \nu_1(\Om^\#),$
 \item if equality happens then $\Om$ is a $\mu$-translate of $\Om^\#$.
\end{enumerate}
\end{theorem}
This theorem generalises the result of Payne-Weinberger \cite{P-W}, and proves that `{\it among the class of domains with given measure, free along the interior boundaries and fixed along a sphere of given radius as the outer boundary, the concentric annulus $\Om^\#$ has the highest fundamental frequency}'. For example, see below  Figure \ref{figure}: 
\begin{figure}[h]\label{figure}
    \centering
    \includegraphics[width=0.75\textwidth]{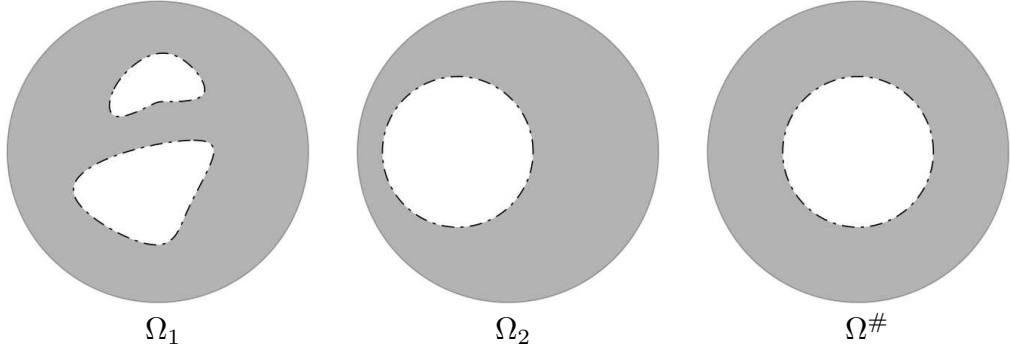}
    \caption{The domains $\Om_1, \Om_2$ and $\Om^\#$ have the same measure  and the same outer boundary (up to a translation). Then  $\nu_1(\Om^\#)$ is strictly greater than $\nu_1(\Om_i).$}
    \end{figure}

       \noi 
       In general, the domain can have any finite number of  interior boundaries. Comparing the values of $\nu_1(\Om_i)$ with the number of interior boundaries is also an interesting problem. In a subsequent article, using the shape derivative  technique,  we also prove that the  $\nu_1(B_1\setminus B_0)$ strictly decreases when the inner ball $B_0$ moves towards the boundary  of  the outer  ball $B_1.$  

 We give two proofs for the above theorem. The first one mainly use some ideas from  \cite{P-W} and an isoperimetric inequality  of the following form:
\begin{equation}
s(\de)^{\pr N}\leq |\Ga_0|^{\pr N} - C(N)v(\de).
\end{equation} Even though Nagy's inequality is not available for higher dimensions, our assumption on  $\Ga_1$  allows us to derive the above isoperimetric inequality. In the second proof, we borrow some ideas from \cite{Hersch}  and give a proof without using the isoperimetric inequality. We feel that the proofs of Payne-Weinberger and Hersch  need more appreciation and may be applicable for other related problems; this is the reason for presenting two proofs for Theorem \ref{thm 1}.

Next we state the  R-F-K inequality for $\tau_1$: 
\begin{theorem}\label{thm 2}
Let $p\in(1,\infty)$  and $\Om$ be as in Theorem \ref{thm 1}. If the inner boundary $\Ga_1$ of $\Om$ be a sphere, then
\begin{enumerate}[(i)]
 \item $ \tau_1(\Om) \le  \tau_1(\Om_\#),$
 \item if equality happens then $\Om$ is a $\mu$-translate of $\Om_\#$.
\end{enumerate}
\end{theorem}
For proving the above result, we introduce a parametrization $t$, analogous to the one in \cite{Hersch}  as
\[t(\de)=\int_0^\de \displaystyle\frac{\dr}{s(r)^{\pr{p}-1}}.\] This parametrization and the eigenfunction corresponding to $\tau_1(\Om_\#)$ helps us to construct a test function on $\Om$, whose Rayleigh quotient is smaller than $\tau_1(\Om_\#).$ 

For $N=2$, as  we have Nagy's inequality for a general multiply connected bounded domain, we extend the results of Payne-Weinberger and Hersch for the $p$-Laplacian with $p\in(1,\infty).$ 
\begin{theorem}\label{thm 3}
Let $p\in(1,\infty)$ and $\Om\subset \R^2$  be a $\C^1$-smooth multiply connected bounded domain.  Then
  \[(i)\ \nu_1(\Om) \le  \nu_1(\Om^\#)\quad \text{and}\quad  (ii)\ \tau_1(\Om) \le  \tau_1(\Om_\#).\]
\end{theorem}
\textbf{The second eigenvalues:}
In \cite{Anoop-Drabek}, Anoop et al. proved that the eigenfunctions corresponding to the second Dirichlet eigenvalue of $p$-Laplacian are non-radial. As an application of above theorems, we prove the similar results for the second eigenvalues of $p$-Laplacian on a ball or an annulus with various boundary conditions. 
Let $\mu_2(\Om)$ be the second (first non-zero) Neumann eigenvalue of $p$-Laplacian; let $\nu_2(\Om)$ and $\tau_2(\Om)$ be the second eigenvalues of \eqref{ND} and \eqref{DN} respectively. 
\noi Now, we have the following results.
\begin{theorem}\label{thm 4}
	For $N\geq 2$ and $1<p<\infty$, let $\Om\subset\RN$ be a ball or a concentric annulus in $\RN$. Then the nodal set of any eigenfunction associated with $\mu_2(\Om)$ cannot be a concentric sphere.
\end{theorem}

\begin{theorem}\label{thm 5}
		For $N\geq 2$ and $1<p<\infty$, let $\Om\subset\RN$ be a concentric annulus. Then the nodal set of any eigenfunction associated with $\nu_2(\Om)$ and $\tau_2(\Om)$ cannot be a concentric sphere. 
\end{theorem}

This article is organised as follows. In the first section, we  briefly illustrate  the interior parallels  and derive Payne-Weinberger's type isoperimetric inequality. The proofs for Theorem \ref{thm 1}, Theorem \ref{thm 2} and Theorem \ref{thm 3} are also given in this section. In the next section, we give the proofs of Theorem \ref{thm 4} and Theorem \ref{thm 5}. In Section \ref{section 4}, we  remark on a similar result for the elasticity problem and state a few open problems related to the mixed  eigenvalue problems. In the appendix, we discuss some properties of the first eigenfunctions and give the variational characterisation of the second eigenvalues.
\section{Proofs of main results}\label{section 2}
In this section we give proofs for  all of our theorems. First, we briefly describe the interior parallels, and then obtain the inequalities analogous to Nagy's inequalities.

\noi{\bf The parallel sets:} For $N\geq 2$, let $\Om$ be a $\C^1$-smooth multiply connected bounded domain in $\R^N$ with interior holes, and let  $\Ga_0$ be its outer boundary and $\Ga_1$  be an inner boundary. For $\de\ge 0$,  consider the following sets: 
 \begin{enumerate}[align=left]
 	\item[$\mathcal{A}_0(\de) =$]  the  set of all points in $\Om$ that are at a distance less than $\de$ from $\Ga_0,$
 	\item[$\mathcal{A}_1(\de)=$] the  set of all points in $\Om$ that are at a distance less than $\de$ from $\Ga_1,$
 	\item[$\mathcal{L}_0(\de)=$] $\pa \mathcal{A}_0(\de)\setminus \Ga_0; \quad \mathcal{L}_1(\de) =\pa \mathcal{A}_1(\de)\setminus \Ga_1.$
 \end{enumerate}
The set $\mathcal{L}_0(\de)$ is known as  an {\it{inner parallel surface}} to $\Ga_0$ and  $\mathcal{L}_1(\de)$ is an {\it{outer parallel surface}} to $\Ga_1$. The sets $\mathcal{L}_0(\de)\cap \Om$ and  $\mathcal{L}_1(\de)\cap\Om$ are collectively called as the {\it interior parallels} in $\Om.$ Since we will be dealing with only one type of parallel surfaces at a time, we use the same notations for the measures of interior parallels:
\[S(\de)=|\mathcal{L}_{i}(\de)| ;\quad s(\de)= |\mathcal{L}_{i}(\de)\cap \Om|,\mbox{ for } i=0,1.\]
    Thus, for $\Om\subset \R^2$ Nagy's inequality gives:
\begin{enumerate}[align=right,]
	\item[] Interior parallels to $\Ga_0$: $\displaystyle s(\de)  \le |\Ga_0|-2\pi \de \mbox{ a.e. on} \left[0,\frac{|\Ga_0|}{2\pi}\right];$
	\item[] Interior parallels to $\Ga_1$: $\displaystyle s(\de) \le |\Ga_1|+2\pi \de \mbox{ a.e. on } [0,\infty).$
\end{enumerate}
\subsection*{\bf The mixed eigenvalue problem with Dirichlet condition on $\Ga_0$:} For this case, we consider the  interior parallels to $\Ga_0$. We make the assumption that $\Ga_0$ is a sphere. Let $|\Ga_0|=N\om_N R_1^{N-1}$ and $\Om^\#=B_{R_1}(0)\setminus\overline{B_{R_0}(0)}$ for some $R_1,R_0>0$. Now the following Nagy's type inequality is immediate: 
\begin{equation}\label{Nagy-out}
S(\de)= N\om_N(R_1-\de)^{N-1}, \de\in [0,R_1].
\end{equation}
Next we derive an isoperimetric inequality in higher dimension, analogous to that of Payne-Weinberger. First, define
\begin{align*}
  \de_\Om = \sup \{\de:&  \mathcal{L}_0(\de)\cap \Om\ne\emptyset\},\\
  v(\de) :=\int_0^\de s(t)\dt&;\quad V(\de):=\int_0^\de S(t)\dt.
 \end{align*}
 Observe that $v(\de)=|\mathcal{A}_0(\de)\cap \Om|,\ V(\de)=|\mathcal{A}_0(\de)|,\ v(\de_\Om)=|\Om|$ and $\de_{\Om^\#}=R_1-R_0$. Also the map $v$ is differentiable a.e., and $\pr v(\de)=s(\de)$ a.e. on $[0,\de_\Om].$ Recall that, for $N\ge 2$, $N^{\prime}=\frac{N}{N-1}$ and $C(N)= N^{N^{\prime}} \om_N^{N^{\prime}-1}$.
\begin{lemma}\label{isoperimetric}
If $\Ga_0$ is a sphere then
\begin{equation}\label{eqn:iso}
 s(\de)^{\pr N}\leq |\Ga_0|^{\pr N}-C(N)v(\de),\, \de\in [0,\de_{\Om}].
\end{equation}
  \end{lemma}
\noi{\it Proof}
As $s(\de)\le S(\de),$ from \eqref{Nagy-out} we have
\begin{align*}
s(\de)&\le N\om_N(R_1-\de)^{N-1},\, \de\in [0,\de_\Om].
\end{align*}
By integrating the above inequality from $0$ to $\de$ yields
\[v(\de)\le  \om_NR_1^N-\om_N(R_1-\de)^N.\] Now multiply the above inequality by $C(N)$, then simple computations leads to the required inequality as below:
\begin{equation}\label{Ineq2:Int}
 C(N)v(\de)\!\le\! N^{\pr N}\om_N^{\pr N}R_1^N-N^{\pr N}\!\om_N^{\pr N}\!(R_1-\de)^N\! =|\Ga_0|^{\pr N}\!-L(\de)^{\pr N}\!\le |\Ga_0|^{\pr N}-s(\de)^{\pr N}.
\end{equation}
\qed
\begin{remark}\label{rmk:isoperimetric}
 If $\Om$ is a $\mu$-translate of $\Om^\#$, then $s(\de)=S(\de) \mbox{ for } \de\in[R_1,R_0]$ and hence the equality holds in the isoperimetric inequality. i.e.,
\begin{equation*}
	s(\de)^{\pr N}= |\Ga_0|^{\pr N}-C(N)v(\de),\, \de\in [0,\de_\Om].
\end{equation*}
\end{remark}
\noi The next lemma will show that, if $\Ga_0$ is a sphere and the equality holds in \eqref{eqn:iso}, then $\Om$ is necessarily a $\mu$-translate of $\Om^\#$.
\begin{lemma}\label{interior_DN} If $\Om^\#$ is not a $\mu$-translate of $\Om$, then
\begin{enumerate}[(i)]
\item  $ R_1-R_0<\de_{\Om},$
\item there exists $\de_0\in [0,\de_\Om]$ such that  $s(\de)^{\pr N}< |\Ga_0|^{\pr N}-C(N)v(\de)  \text{ for } \de\ge \de_0.$
\end{enumerate}
\end{lemma}
\noi{\it Proof.}
	$(i)$ Since the outer boundary of $\Om$ is a sphere without loss of generality, we may assume that the outer boundaries of  $\Om$ and $\Om^\#$ are the same. If $R_1-R_0\ge \de_\Om$ then $\Om$ must be a subset of $B_{R_1}(0)\setminus \overline{B_{R_1-\de_\Om}(0)}\subseteq \Om^\#,$ since  $\mathcal{L}_0(\de)\cap \Om= \emptyset$ for $\de>\de_\Om$. If the strict inequality holds then the  inclusion is  strict, and if the equality holds, then  $\Om$ must be a $\mu$-translate of $\Om^\#.$ In both cases, we get a contradiction and hence $R_1-R_0< \de_\Om.$\\
	$(ii)$ If $v = V$ on $[0,\de_\Om],$ then $|\Om|=V(R_1-R_0)<V(\de_\Om)=v(\de_\Om)=|\Om|$. A contradiction and hence there  exists $\de_0\in [0,\de_\Om]$ such that $v(\de_0)< V(\de_0).$ Since $s(\de)\le S(\de)$, we get $v(\de)< V(\de),$ for $\de>\de_0.$
	Now the same calculations as in \eqref{Ineq2:Int}, gives the required strict inequality.
	
	\qed

\noi From Remark \ref{rmk:isoperimetric}, we have  
$ s(\de)^{\pr{N} }\le |\Ga_0|^{\pr{N}}-C(N)v(\de), \, \de\in[0,\de_{\Om}].$
Motivated by this inequality, we define a function $r$ as below:
\begin{equation}\label{transform}
C(N)\om_N r(\de)^{N }= |\Ga_0|^{\pr{N}}-C(N)v(\de),\, \de\in[0,\de_{\Om}].
\end{equation}
Observe that,  $r(0)=R_1,\, r(\de_{\Om})=R_0$ and $r$ is strictly decreasing on $[0,\de_{\Om}]$. Thus $r$ maps $[0,\de_\Om]$ onto $[R_0,R_1]$.
\begin{lemma}\label{second proof}
Let $r$ be defined as above. Then the map $r$ is differentiable and   
$|\pr r(\de)| \leq 1$.
\end{lemma}
\noi{\it Proof}
By differentiating \eqref{transform} we get $N C(N)\om_N r(\de)^{N-1}\pr r(\de)= -C(N)s(\de).$ Therefore,
\[|\pr r(\de)|^{\pr N}=\frac{s(\de)^{\pr{N}}}{N^{\pr N}\om_N^{\pr N}r(\de)^N}=\frac{s(\de)^{\pr{N}}}{|\Ga_0|^{\pr{N}}-C(N)v(\de)}.\]
Now the conclusion follows easily from \eqref{eqn:iso}.

\qed

\noi Next we  give our first proof for Theorem \ref{thm 1}, along the same lines in the proof of Payne-Weinberger (Section-II of \cite{P-W}).

\noi{\bf{The first proof of Theorem \ref{thm 1}}.}\\
$(i)$ Let $u_0\in W^{1,p}(\Om^\#)$ be an eigenfunction corresponding to $\nu_1(\Om^\#)$. Since $u_0$ is radial, there exits $\phi\in C^1(\R)$ such that $\phi(R_1)=0$ and $u_0(x)=\phi(|x|)$. Now consider the  function  $u(x)= \phi(r(\de(x))),$ for $x\in \Om$. Observe that  $u\in W^{1,p}(\Om)$  and $u|_{\Ga_0}=0.$  Further,
\[|\Gr u(x)|=\left|\pr \phi(r(\de(x)))\right||\pr r(\de)| |\Gr \de|\leq |\pr \phi(r(\de(x)))|, \ (\mbox{since } |\pr r(\de)|\leq 1,\, |\nabla \de|=1).\]
Therefore,
\begin{eqnarray*}
\int_{\Om}|\Gr u(x)|^p\ \! \dx&\leq & \int_{\Om}\left|\pr \phi(r(\de(x)))\right|^p \dx\\ &=& \int_{0}^{\de_\Om}\Big(\int\limits_{{\mathcal{L}}(\de)}\left|\pr \phi(r(\de(x)))\right|^p \dsig \Big) \dd = \int_{0}^{\de_\Om}\left|\pr \phi(r(\de))\right|^p s(\de) \dd.
\end{eqnarray*}
From \eqref{transform} we have $N\om_N r^{N-1}\dr = - s(\de)\dd.$ Now by setting $r(\de)=r$ we get
\[\int_{\Om}|\Gr u(x)|^p\ \! \dx \leq N\om_N \int_{R_0}^{R_1}\left|\pr \phi(r)\right|^p r^{N-1} \, \dr.\]
Similarly, 
\[\displaystyle\int_{\Om^\#}|\Gr u_0(x)|^p\! \dx=N\om_N \int_{R_0}^{R_1}\left|\pr \phi(r)\right|^p r^{N-1} \ \! \dr.\]
Further, 
\[\displaystyle\int_\Om |u(x)|^p\dx = N\om_N \int_{R_0}^{R_1}|\phi(r)|^p r^{N-1}\dr= \displaystyle\int_{\Om^\#} |u_0(x)|^p\dx.\]
Combining the above inequalities, we conclude that
\[\nu_1(\Om)\leq \nu_1(\Om^\#).\]
$(ii)$ If $\Om$ is not a $\mu$-translate of $\Om^\#$, then by  Lemma \ref{interior_DN} there exists $\de\in (0,\de_\Om)$ such that  
$s(\de)^{N^{\prime}}< |\Ga_0|^{N^{\prime}}-C(N)v(\de), \, \de\in (\de_0,\de_\Om).$ Now from Lemma \ref{second proof}, we have $|\pr r(\de)|< 1$ for $ \de\in (\de_0,\de_\Om).$ Therefore,
\[\displaystyle\int_{\Om}|\Gr u(x)|^p \dx <N\om_N \int_{R_0}^{R_1}\left|\pr \phi(r)\right|^p r^{N-1} \ \! \dr = \int_{\Om^\#}|\Gr u_0(x)| \ \! \dx.\] 
Thus $\nu_1(\Om)< \nu_1(\Om^\#).$

\qed

Next we give another proof of Theorem \ref{thm 1}, without (explicitly) using the isoperimetric inequality \eqref{eqn:iso}. For this, first  we prove the following lemma.
\begin{lemma}\label{lemma:outer}  Let $h(\al)= s(v^{-1}(\al))$ and $H(\al)=S(V^{-1}(\al))$ for $\al\in[0,|\Om|].$ Then 
\begin{enumerate}[(i)]
\item  $ h(\al)\leq H(\al), \ \al\in[0,|\Om|],$
\item if $\Om$ is not a $\mu$-translate of $\Om^\#$, then there exists $\al_0\in [0,|\Om|]$ such that $h(\al)< H(\al),\ \al\geq \al_0.$
\end{enumerate}
\end{lemma}
\noi{\it Proof}
$(i)$ Since $v(\de)\leq V(\de)$ and $S$ is monotonically decreasing, we have $v^{-1}(\al)\geq V^{-1}(\al)$ and hence $S(v^{-1}(\al))\leq S(V^{-1}(\al)).$ As $s(\de)\le S(\de)$ we get the required conclusion.\\
$(ii)$ If $\Om$ is not a $\mu$-translate of $\Om^\#$,  then as in the proof of $(ii)$ of Lemma \ref{interior_DN}, there exists a $\de_0\in [0,\de_\Om]$ such that $v(\de) < V(\de)$ for all $\de\geq \de_0.$ By setting $\al_0=V(\de_0)$, we obtain the strict inequality  as $S$ is strictly decreasing.

\qed

\noi{\bf The second proof of Theorem \ref{thm 1}.}\\
$(i)$ Let $u_0\in W^{1,p}(\Om^\#)$ be the eigenfunction corresponding to $\nu_1(\Om^\#)$ with $u_0>0$. Since $u_0$ is radial, there exists a $\psi\in C^1(\R)$ such that $\psi(0)=0$ and  
\[u_0(x)=\psi(R_1-|x|)=(\psi\circ V^{-1})\left(V(\de(x))\right),\ x\in \Om^\#.\]
By setting $\phi=\psi\circ V^{-1},$ we have $u_0(x)=\phi\left(V(\de(x))\right).$ Clearly, $\phi\in C^1(\R)$ with $\phi(0)=0.$ Now, define a function $u$ on $\Om$ as $u(x)=\phi(v(\de(x))),\ x\in\Om.$ Observe that, $u\in W^{1,p}(\Om), \, u|_{\Ga_0}=0,$ { and } 

\[|\Gr u(x)|=|\pr\phi\left(v(\de(x))\right)\!|\, s (\de(x)), \mbox{a.e. on } \Om.\]  Therefore,
 \[\int_{\Om}|\Gr u(x)|^p \dx = \int_0^{\de_{\Om}} \left(\int_{\mathcal{L}_{\de}} |\Gr u(x)|^p\dsig \right) \, \dd= \int_0^{\de_\Om} |\pr\phi\left(v(\de)\right)|^p s(\de)^{p+1} \dd.\]
Using the change of variable $\al=v(\de)$, we deduce
 \[\int_{\Om}|\Gr u(x)|^p \dx = \int_0^{|\Om|} |\pr\phi(\al)|^p h(\al)^p \da.\]
A similar calculation yields, 
\begin{align*}
                      \int_{\Om^\#}|\Gr u_0(x)|^p \dx =\int_0^{|\Om|}|\pr\phi(\al)|^p H(\al)^p \da.
                     \end{align*}
Further, $\displaystyle\int_{\Om^\#}| u_0(x)|^p \dx=\int_0^{|\Om|}\phi(\al)^p\da= \int_{\Om}| u(x)|^p \dx.$ By Lemma \ref{lemma:outer}, $h(\al)\leq H(\al)$, so we get 
\[\displaystyle\int_0^{|\Om|}\!\! |\pr\phi(\al)|^p  h(\al)^p \da \leq \int_0^{|\Om|}\!|\pr\phi(\al)|^p H(\al)^p \da.\]
Now, by the variational characterisation of $\nu_1(\Om)$, we conclude $\nu_1(\Om)\leq \nu_1(\Om^\#)$.

\noi$(ii)$ If $\Om$ is not a $\mu$-translate of $\Om^\#$, then $(ii)$ of Lemma \ref{lemma:outer} and the above calculations gives $\nu_1(\Om)< \nu_1(\Om^\#)$.
\qed
\begin{remark}
 {\rm Notice that, in our case ($\Ga_0$ is a sphere), the equality in the isoperimetric inequality implies  $\Om$ is a $\mu$-translate of $\Om^\#$ and hence $\nu_1(\Om)=\nu_1(\Om^\#).$}
\end{remark}
\subsection{\bf The mixed eigenvalue problem with Dirichlet condition on $\Ga_1$:}
In this case, we consider the interior parallel surfaces from the inner boundary $\Ga_1$. Let $\de_\Om:= \sup \{\de:  \mathcal{L}_1(\de)\cap \Om\ne\emptyset\}$. We assume that $\Ga_1$ is a sphere. Let $|\Ga_1|=N\om_N R_0^{N-1}$ and $\Om_\#=B_{R_1}(0)\setminus\overline{B_{R_0}(0)}$ for some $R_1,R_0>0$. Now, we clearly have the following Nagy's type inequality:
\begin{equation}\label{Nagy-in}
S(\de)=N\om_N\left (R_0+\de\right)^{N-1}, \ \de\in[0,\de_\Om]. 
\end{equation}
For $p\in(1,\infty)$, we define a parametrization $t$ as 
\[t(\de)=\int_0^\de \frac{\dr}{s(r)^{\pr{p}-1}};\qquad T(\de)=\int_0^\de \frac{\dr}{S(r)^{\pr{p}-1}}, \ \de\in[0,\de_\Om],\]
where $\pr{p}=\frac{p}{p-1}$ is the H\"{o}lder conjugate of $p$. Let $T_\#:=T(R_1-R_0);\; t_\Om:= t(\de_\Om)$. Notice that $t,T$ and $S$ are strictly increasing and $ s(\de)\le S(\de).$ Now consider the maps $g(\al)=s(t^{-1}(\al)) \mbox{ and } G(\al)=S(T^{-1}(\al)) \mbox{ on } [0,|\Om|].$
\begin{lemma}\label{Inner1}
Let $p\in(1,\infty)$ and $\Om$ be as in Theorem \ref{thm 2}.  If $\Ga_1$ is a sphere, then 
 \begin{enumerate}[(i)]
 \item $\de_\Om\ge R_1-R_0$ and $T_\#\le t_\Om,$
 \item $g(\al)\le G(\al),\ \al\, \in (0,T_\#),$
 \item $\displaystyle\int_0^{t_\Om}g(\al)^{\pr{p}}\da=\int_0^{T_\#}G(\al)^{\pr{p}}\da=|\Om|$.
   \end{enumerate}
\end{lemma}
\noi{\it Proof} $(i)$ If $\de_\Om< R_1-R_0,$ then $\Om$ must be a $\mu$-translate of a subset of the annulus  $B_{R_0+\de_\Om}(0)\setminus \overline{B_{R_0}(0)}.$ A contradiction and hence $\de_\Om\ge R_1-R_0.$ Now as $T(\de)\le t(\de),$ we get 
$T_\#\le t_\Om.$\\
$(ii)$ For $\al\in [0,T_\#]$, clearly $t^{-1}(\al)\le T^{-1}(\al)$ and hence \[g(\al)=s(t^{-1}(\al))\le S(t^{-1}(\al))\le S(T^{-1}(\al))=G(\al).\]
$(iii)$ By changing the variable $\al=t(\de)$, we get 
\[\displaystyle \int_0^{t_\Om} g(\al)\da=\int_0^{t_\Om}s(t^{-1}(\al))^{\pr p} \da= \int_0^{\de_\Om}s(\de)^{\pr p}\frac{1}{s(\de)^{\pr p -1}}\dd = \int_0^{\de_\Om}s(\de)\dd=|\Om|.\] Similarly $\displaystyle\int_0^{T_\#}G(\al)^{\pr p}\da=|\Om_\#|=|\Om|.$
\qed   
\begin{lemma}\label{lemma:inner} Let $\Om$ be as in Theorem \ref{thm 2}.
 If $\Om$ is not a $\mu$ translate of $\Om_\#$, then 
 \begin{enumerate}[(i)]
  \item $\de_\Om>R_1-R_0,$
  \item $T_\#< t_\Om,$
  \item there exists $\al_0\subset [0,T_\#]$ such that $g(\al)< G(\al), \ \al \in (\al_0,T_\#].$ 
 \end{enumerate}
\end{lemma}
\noi{\it {Proof}}
$(i)$ By the previous lemma, $\de_\Om\geq R_1-R_0$. If $\de_\Om= R_1-R_0,$ then $\Om$ must be a  $\mu$-translate of the annulus  $\Om_\#.$ Hence $\de_\Om> R_1-R_0.$\\
$(ii)$ Follows from $(i)$ as $T$ and $t$ are strictly increasing.\\
$(iii)$ If $t=T$ on $[0,T_\#]$, then $s=S$ a.e. on $[0,T_\#]$ and this implies that $|\Om|=\displaystyle\int_0^{R_1-R_0}S(r)\dr \leq \int_0^{R_1-R_0}s(r)\dr <\int_0^{\de_\Om}s(r)\dr=|\Om|$. A contradiction and hence there exists a $\de_0\in[0,\de_\Om]$ such that $t(\de_0)>T(\de_0).$ Since $s(\de)\leq S(\de)$, we get $t(\de)>T(\de)$, for $\de>\de_0$. By setting $\al_0=v(\de_0)$, we obtain the required result.

\qed
\noi {\bf Proof of Theorem \ref{thm 2}.}\\
$(i)$ Let $u_1$ be a non-negative eigenfunction corresponding to $\tau_1(\Om_\#)$. Since $u_1$ is radial in $\Om_\#$, there exists $\psi\in C^1(\R)$ such that $\psi(0)=0$ and 
\[u_1(x)=\psi(|x|-R_0)=\psi\circ T^{-1}\big(T(\de(x)\big)=\phi(T(\de(x))), \ x\in\Om_\#,\]
where $\phi=\psi\circ T^{-1}$. Notice that $u_1$ has its maximum on the outer boundary $\Ga_0$, and hence
\begin{equation}\label{inner:eqn1}                        \phi(T_\#)\ge \phi(\al), \, \al\in [0, T_\#].
\end{equation}
Now, define a function $u$ on $\Om$ as 
\begin{eqnarray*}
 u(x)=\left \{ \begin{array}{ll} \phi\big(t(\de(x))\big),\; t(\de(x))\in[0,T_\#],\\
        \phi(T_\#),\; t(\de(x))\in(T_\#,t_\Om].
        \end{array} \right.
\end{eqnarray*}
Clearly, $u\in W^{1,p}(\Om)$ and $u|_{\Ga_1}=0$. Further,
\[\Gr u(x)=\left \{ \begin{array}{ll}\displaystyle \frac{\phi^{\prime}(t(\de(x)))}{s(\de(x))^{\pr{p}-1}}\Gr \de(x),&  t(\de(x))\in[0,T_\#],\\
        0,& t(\de(x))\in(T_\#,t_\Om].
        \end{array} \right. \]
Therefore,
\begin{eqnarray*}
  \int_{\Om}|\Gr u(x)|^p\dx &=& \displaystyle\int_0^{t^{-1}(T_\#)}\left(\int_{{\mathcal{L}}(\de)}|\Gr u(t(\de(x))|^p\dsig\right) \dd\\ &=& \int_0^{t^{-1}(T_\#)}|\pr{\phi}\big(t(\de)\big)|^p \frac{\dd}{s(\de)^{\pr{p}-1}}.
 \end{eqnarray*}
  By the change of variable $\al=t(\de)$, we get 
  $\displaystyle\int_{\Om}|\Gr u(x)|^p\dx=\displaystyle\int_0^{{T_\#}}|\pr{\phi}(\al)|^p \da.$\\
Further, we obtain 
\begin{eqnarray*}
\int_{\Om_\#}|\Gr u_1(x)|^p\dx &=&\int_0^{{T_\#}}|\pr{\phi}(\al)|^p \da\,, \int_{\Om_\#}|u_1(x)|^p\dx = \int_0^{{T_\#}}|\phi(\al)|^p G(\al)^{\pr{p}} \da,\\
          \int_{\Om}|u(x)|^p\dx & =& \int_0^{T_\#}|\phi(\al)|^p g(\al)^{\pr{p}} \da+|\phi(T_\#)|^p \int_{T_\#}^{t_\Om} g(\al)^{\pr{p}} \da.          
         \end{eqnarray*}
Using the facts $g(\al)\leq G(\al)$ and $\phi(t)\leq \phi(T_\#)$ (by Lemma \ref{Inner1} and \eqref{inner:eqn1}), we estimate the following:
\begin{align*}
              \int_{\Om_\#}|u_1(x)|^p\dx&-\int_{\Om}|u(x)|^p\dx\\ &=\displaystyle\int_0^{T_\#}\!|\phi(\al)|^p \Big(G(\al)^{\pr{p}}-g(\al)^{\pr{p}}\Big) \da - \int_{T_\#}^{t_\Om}|\phi(T_\#)|^p g(\al)^{\pr{p}}\\
                            &\le \displaystyle|\phi(T_\#)|^p\left( \int_0^{{T_\#}}\Big(G(\al)^{\pr{p}}-g(\al)^{\pr{p}}\Big) \da - \int_{T_\#}^{t_\Om} g(\al)^{\pr{p}}\right)=0.
             \end{align*}
             Therefore, we obtain $\displaystyle \int_{\Om_\#}|u_1(x)|^p\dx \le \int_{\Om}|u(x)|^p\dx.$ Now  the above inequalities and  the variational characterisation of $\tau_1$ gives 
             \[\tau_1(\Om)\leq \tau_1(\Om_\#).\]
$(ii)$ If $\Om$ is not a $\mu$-translate of $\Om_\#$, then from Lemma \ref{lemma:inner}, there exists $\al_0\in [0,T_\#]$ such that  $g(\al)<G(\al), \, \al \in (\al_0,T_\#].$ Therefore, $\displaystyle\int_{\Om_\#}|u_1(x)|^p\dx<\int_{\Om}|u(x)|^p\dx$ and hence $\tau_1(\Om)<\tau_1(\Om_\#).$

\qed
\begin{remark} \rm{ Consider the annular domains of the form $B\setminus \overline{B_0}\subset \R^N$ with $N\ge 2.$ For $p\in (1,\infty)$, the first eigenvalues $\nu_1(B\setminus \overline{B_0})$ and $\tau_1(B\setminus \overline{B_0})$ are maximum only if the balls are concentric.}
\end{remark}
\subsection{The mixed eigenvalue problems in dimension 2:}
 For both the annular regions $\Om^\#$ and $\Om_\#,$ we denote the lengths of the interior parallels by $\overline{S}(\de)$.

\noi{\bf Proof of Theorem \ref{thm 3}.}\\
$(i)$ In this case, we consider the interior parallels from $\Ga_0.$ Thus we have $\overline{S}(\de)=|\Ga_0|-2\pi\de$. Further, we define 
\[ \overline{V}(\de)=\int_0^\de \overline{S}(t) \dt,\, \de\in[0,R_1-R_0];\quad \overline{H}(\al)=\overline{S}\big(\overline{V}^{-1}(\al)\big),\, \al\in [0,|\Om|],\]
 As before, we obtain $ \overline{S}(\de)^2=|\Ga_0|^2-4\pi\overline{V}(\de),\, \de\in [0,R_1-R_0]$ and this yields
\[\overline{H}(\al)^2=|\Ga_0|^2-4\pi\al,\, \al\in [0,|\Om|].\]
From  Nagy's inequality, we have $s(\de)\le \overline{S}(\de)$ and this gives the following isoperimetric inequality
\[s(\de)^2\leq |\Ga_0|^2-4\pi v(\de), \, \de\in [0, \de_\Om].\] Therefore, $h(\al)^2\leq \overline{H}(\al)^2,\, \al \in [0,|\Om|].$  Now the proof follows using a similar set of arguments as in the second proof of part $(i)$ of Theorem \ref{thm 1}.\\
$(ii)$ For this case, we consider the interior parallels from $\Ga_1.$ Thus we have $\overline{S}(\de)=|\Ga_1|+2\pi\de, \, \de\in (0,\de_\Om)$ and Nagy's inequality gives $s(\de)\le \overline{S}(\de).$ Now define a parametrization as
\[\overline{T}(\de):=\int_0^\de \frac{\dt}{\overline{S}(t)^{\pr p -1}}, \, \de\in[0,\de_\Om].\]
Let $\overline{G}(\al)=\overline{S}\big(\overline{T}^{-1}(\al)\big),\, \al\in [0,T_\#].$ Now from $(ii)$ of Lemma \ref{Inner1}, we get $g(\al)\leq \overline{G}(\al),\, \al\in [0,T_\#].$ Now the rest of the proof is same as the proof of part $(i)$ of Theorem \ref{thm 2}.

\qed
\section{Applications}\label{section 3}
In this section we prove Theorem \ref{thm 4} and Theorem \ref{thm 5}.
First, we give a variational characterisation of the second eigenvalue. For this, let
\[\M\dq \S\cap\left\{u\in\wp :\int_{\Om}|u|^{p-2}u=0 \right\},\]
\[\F_2\dq\left\{h(S^{1}): h \mbox{ is an odd continuous map from }S^{1} \mbox{ into }\S\cap W_{\Ga_0} \right\},\]
\[\G_2\dq\left\{h(S^{1}): h \mbox{ is an odd continuous map from }S^{1} \mbox{ into }\S\cap W_{\Ga_1} \right\}.\]
Then 
\[\mu_2(\Om)=\inf\limits_{u\in\M} J(u),\ \nu_2(\Om)=\inf\limits_{\A\in\F_2}\sup\limits_{u\in\A} J(u),\ \tau_2(\Om)=\inf\limits_{\A\in\G_2}\sup\limits_{u\in\A} J(u). \]
\noi{\bf  Proof of Theorem \ref{thm 4}.} By the translation invariance of $p$-Laplacian, we can take $\Om=B_{R_1}(0)\setminus \overline{B_{R_0}(0)}$ for some $0\leq R_0<R_1<\infty$. 
 Let $u$ be the eigenfunction associated with $\mu_2(\Om)$, then $u$ must change its sign in $\Om$. Suppose the nodal set $\{x\in\Om: u(x)=0\}$ of $u$ is a sphere of radius $r\in(R_0,R_1)$ centred at origin. Then we have
 \[\tau_1(B_r(0)\setminus\overline{B_{R_0}(0)})=\mu_2(\Om);\qquad \nu_1(B_{R_1}(0)\setminus\overline{B_{r}(0)})=\mu_2(\Om).\]
Now for $s\in (0,r-R_0)$, by Theorem \ref{thm 1} and Theorem \ref{thm 2} we have
\[\tau_1(B_r(se_1)\setminus\overline{B_{R_0}(0)})\leq \mu_2(\Om)\mbox{ and } \nu_1(B_{R_1}(0)\setminus\overline{B_r(se_1)})<\mu_2(\Om).\] 
 Let $\phi_1$ and $ \phi_2$ be positive eigenfunctions corresponding to $\tau_1(B_r(se_1)\setminus\overline{B_{R_0}(0)})$ and $\nu_1(B_{R_1}(0)\setminus\overline{B_r(se_1)})$ respectively, with the normalisation $\int|\phi_1|^{p-1}=1=\int|\phi_2|^{p-1}$. Let $\tilde{\phi}$ denotes the extension of $\phi$ to $\Om$ by zero. Then $\phi=\tilde{\phi_1}-\tilde{\phi_2}\in \wp$ with $\int_{\Om}|\phi|^{p-2}\phi=0$, and
	\begin{align*}
	J(\phi)=\int_{\Om}|\nabla \phi|^p&=\int|\nabla \phi_1|^p+\int|\nabla \phi_2|^p\\
	&= \tau_1(B_r(se_1)\setminus\overline{B_{R_0}(0)})\int|\phi_1|^p+\nu_1(B_{R_1}(0)\setminus\overline{B_r(se_1)})\int|\phi_2|^p\\
	&<\mu_2(\Om)\left(\int|\phi_1|^p+\int|\phi_2|^p\right)=\mu_2(\Om)\int_{\Om}|\phi|^p.
	\end{align*}
This contradicts the variational characterisation of $\mu_2(\Om)$.

\noindent {\bf Proof of Theorem \ref{thm 5}.} We give the result only for $\nu_2(\Om)$, a similar proof holds for $\tau_2(\Om)$. By the translation invariance of $p$-Laplacian, we can take $\Om=B_{R_1}(0)\setminus \overline{B_{R_0}(0)}$ for some $0< R_0<R_1<\infty$.  Let $u$ be an eigenfunction associated with $\nu_2(\Om)$.
	Suppose the nodal set $\{x\in\Om:u(x)=0\}$ of $u$ is a sphere of radius $r\in(R_0,R_1)$ centred at the origin. Then we have 
 	
	\[\la_1(B_r(0)\setminus\overline{B_{R_0}(0)})=\nu_2(\Om);\qquad  \tau_1(B_{R_1}(0)\setminus\overline{B_{r}(0)})=\nu_2(\Om).\] Thus, for $s\in (0,r-R_0)$,  by Theorem 1.1 of \cite{Anoop} and Theorem \ref{thm 2} we have 
	
	\[\la_1(B_r(se_1)\setminus\overline{B_{R_0}(0)})< \nu_2(\Om) \;\mbox{ and }\; \tau_1(B_{R_1}(0)\setminus\overline{B_r(se_1)})<\nu_2(\Om).\] 
	Let $\phi_1$ and $ \phi_2$ be the positive eigenfunctions corresponding to the eigenvalues $\la_1(B_r(se_1)\setminus\overline{B_{R_0}(0)})$ and $\tau_1(B_{R_1}(0)\setminus\overline{B_r(se_1)})$ respectively, with the normalisation $\int|\phi_1|^{p}=1=\int|\phi_2|^{p}$. Now consider the set $\A\dq \{a\tilde{\phi_1}+b\tilde{\phi_2}: |a|^p+|b|^p=1\}$. Then $\A\in\F_2$ and for $\phi\in\A$ we have
	\begin{align*}
	J(\phi)=\int_{\Om}|\nabla \phi|^p&=\int|\nabla \phi_1|^p+\int|\nabla \phi_2|^p\\
	&= \la_1(B_r(se_1)\setminus\overline{B_{R_0}(0)})|a|^p+\tau_1(B_{R_1}(0)\setminus\overline{B_r(se_1)})|b|^p\\
	&<\nu_2(\Om)\left(|a|^p+|b|^p\right)=\nu_2(\Om).
	\end{align*}
A contradiction to the variational characterisation of $\nu_2(\Om)$. 
\qed
\section{Some remarks and open problems} \label{section 4}
\begin{remark}\label{Remark-Steiner}
 \rm{ Let $\Om$ be a doubly connected planar domain whose inner hole is a convex set and the outer boundary is a parallel curve to the inner boundary. Now one can use Steiner's formula for the convex domain and get the equality in Nagy's inequality. i.e., $S(\de)=|\Ga_1|+2\pi\de.$ For such a domain $\Om,$  without being a $\mu$-translate of $\Om_\#,$ indeed we have  $T_\#=t_\Om$ and $\overline{G}(\al)=g(\al).$ However,  it is not clear whether $\tau_1(\Om)=\tau_1(\Om_\#)$ or not.}
\end{remark}

\begin{remark}[The elasticity problem] For $p\in(1,\infty)$ and for a multiply connected domain $\Om$ in $\R^N$ with $N\ge 2$, one can consider the following elasticity problem:
 
 \begin{equation*}
	\left.
	\hs{-0.8cm}
		\begin{aligned}
	 &-\De_p u= \La |u|^{p-2}u \ \text{in } \, \Om,\\
	  & |\nabla u|^{p-2}\frac{\pa u}{\pa \eta}+k|u|^{p-2}u =0 \text{ on }  \Ga_0,\\
	 &\frac{\pa{u}}{\pa \eta}=0  \text{ on }   \Ga_1\cup\left(\cup_{i=1}^n\ga_i\right),
	\end{aligned}
	\right\}
	\end{equation*}
	where $k>0$ is the elasticity constant. The first eigenvalue $\La_1(\Om)$ has the following variational characterisation:
	
	\[\La_1(\Om)= \displaystyle\inf \left \{\frac{\int_{\Om}\md{\nabla u}^p+ k \int_{\Ga_1}|u|^p}{\int_{\Om}\md{u}^p}: u\in W^{1,p}(\Om)\setminus\{0\} \right\}.\]
	 Then using the similar set of arguments as in the proof of  Theorem \ref{thm 1} one can show that 
	if $\Ga_0$ is  a sphere, then $ \La_1(\Om) \le  \La_1(\Om^\#)$ and the equality holds if and only if $\Om$ is a $\mu$-translate of $\Om^\#$.
\end{remark}

Next, we state a few open problems that are related to the mixed eigenvalue problems that we considered in this article. For more open problems related to the extremum of eigenvalues of various operators, we refer to the book \cite{Henrot}.

\noi \textbf{Open problems:}
\begin{enumerate}[(i)]
 \item {The uniqueness for $N=2$:}  Our results ensures the uniqueness of the domain (up to $\mu$-translates) that gives the equality in the reverse F-K inequality. However, for the general domains ($N=2$), Payne-Weinberger's or Hersch's results are not making any  claim on the uniqueness of the domain. In the best of our knowledge this question is open for $N=2$ and for every $p$. See also our Remark \ref{Remark-Steiner}.
 \item {The mixed  eigenvalue problems for the general multiply connected domains:} We proved our results under the assumptions that the boundary on which the Dirichlet condition is specified are spheres.  We feels that the isoperimetric inequalities may obtained for other domains. For an excellent review on isoperimetric inequalities and related results, we refer to \cite{Osserman,Polya-Szego}. 
 \item {The lower bounds for the first mixed eigenvalues $\nu_1(\Om)$:} Is it possible to find a constant $C$ (depending on the domain) such that  $C \nu_1(\Om^\#)\le \nu_1(\Om)$ and  hence  an upper bound for the isoperimetric deficit $\frac{\nu_1(\Om^\#)-\nu_1(\Om)}{\nu_1(\Om^\#)}.$ The similar problems for Dirichlet eigenvalue are studied in \cite{BCT} and \cite{DFG} for various operators.
 
 \item  {The Dirichlet eigenvalue problem for the general multiply connected domains:} Except for the annular region, the R-F-K  is open for the general doubly connected domains  in the higher dimension, even for the case one of the boundaries is a sphere. 
 
 \item  One can also study the mixed eigenvalue problems with the Dirichlet condition specified on more than one interior boundaries. An upper bound for the case when $N=2$, $p=2$ is given by Hersch, see Section 3 of \cite{Hersch}. 
\end{enumerate}
\appendix
\noi\textbf{Appendices}
\section{The existence of the first eigenvalue and some of its properties}\label{Appendix A}
For $N\geq 2$, let $\Om\subseteq\RN$ be a multiply connected smooth domain such that $\pa \Om =\Ga_0 \sqcup \Ga_1$ with $|\Ga_0|>0$, the $(N-1)$-dimensional measure. For $1<p<\infty$, consider the following mixed eigenvalue problem:
\begin{equation}\label{Genp}
\left.\begin{aligned}
-\De_p u &= \al |u|^{p-2}u  \quad  \text{in }  \Om ,\\
u&=0 \text{ on } \Ga_0, \\
\frac{\pa{u}}{\pa{n}}&=0  \text{ on }  \Ga_1.
\end{aligned}
\right\}
\end{equation}
We say a real number $\al$ is an eigenvalue of \eqref{Genp} if there exists a $u\in W_{\Ga_0}\setminus \{0\}$ such that
\[\int_{\Om}\md{\nabla u}^{p-2}\nabla u\cdot\nabla v = \al \int_{\Om} \md{u}^{p-2}uv,\ \text{ for all } v\in W_{\Ga_0}, \]
where $W_{\Ga_0}= \left\{u\in \wp:\ u|_{\Ga_0}=0 \right\}.$ Recall that
\[J(u)=\int_{\Om}|\nabla u|^p\dx,\ \S=\left\{ u\in \wp: \int_\Om|u|^p\dx=1\right\}.\]
It is easy to verify that the eigenvalues of \eqref{Genp} are precisely the critical values of $J$ on $\S\cap W_{\Ga_0}$.
\begin{proposition}\label{propo:appendix}
	Let $p\in(1,\infty)\mbox{ and } \Om$ as above. Let $\al_1 =\displaystyle\inf_{u\in \S\cap W_{\Ga_0}} J(u).$ Then 
	\begin{enumerate}[$(i)$]
		\item $\al_1$ is an eigenvalue of \eqref{Genp},
		\item $\al_1$ is simple and principal.
	\end{enumerate}
\end{proposition}
{\it Proof.} $(i)$ Let $\left(u_n\right)_{n\in\N}$ be a minimising sequence, i.e., $J(u_n)=\int_{\Om}\md{\nabla u_n}^p\rightarrow \al_1$. 
Then $\left(u_n\right)$ is a bounded sequence in $W_{\Ga_0}$ and hence by the reflexivity, $u_n\rightharpoonup u$ in $W_{\Ga_0}$ for some $u\in W_{\Ga_0}.$ From the compact embedding $W_{\Ga_0}\hookrightarrow L^p(\Om), \ u_n\rightarrow u$ in $L^p(\Om),$ so $\int_{\Om} \md{u}^p=1$ hence $u\in \S.$ Now, by the weakly lower semi continuity of $J$ we have
\[\al_1\leq J(u)\leq \liminf\limits_{n\rightarrow \infty} J(u_n) =\al_1.\]
Hence $\al_1$ is attained for some $u\in W_{\Ga_0}$. Now it is easy to verify that $\al_1$ is a critical value of $J$ and hence an eigenvalue of \eqref{Genp}. If $\al_1=0$ then $\nabla u =0$, this will imply that $u\equiv 0$, as $u\equiv 0$ on the boundary, a contradiction.  Hence $\al_1$ is the first non-zero eigenvalue.\\
$(ii)$ To prove the simplicity of $\al_1$ we use the Picone's identity. Let us recall the Picone's identity: {\it Let $u\geq 0,\ v>0$ a.e., such that $\md{\nabla u},\ \md{\nabla v}$ exists. Define 
	\[L(u,v)\dq \md{\nabla u}^p+(p-1)\displaystyle\frac{u^p}{v^p}\md{\nabla v}^p-p\displaystyle\frac{u^{p-1}}{v^{p-1}}\md{\nabla v}^{p-2}\nabla v,\] and 
	\[R(u,v)\dq \md{\nabla u}^p-\md{\nabla v}^{p-2} \nabla \left(\displaystyle\frac{u^p}{v^{p-1}}\right)\cdot\nabla v.\] Then we have the following: 
	\[L(u,v)=R(u,v) \mbox{ a.e.}\] Furthermore,  $L(u,v)\geq 0$ and $L(u,v)=0$ if and only if $u,v$ are constant multiples on each connected component of $\Om$.}

\noi Let $u,v>0$ be two eigenfunctions corresponding to $\al_1$. By Picone's identity we get
\[0\leq \int_{\Om} L(u,v)\dx=\int_\Om R(u,v)\dx= \displaystyle\int_{\Om}(\md{\nabla u}^p- \al_1\md{u}^p)\dx =0.\] Therefore $u=kv$ for some $k>0$, hence $\al_1$ is simple.

\noi For principality of $\al_1$, let $u\in \S\cap W_{\Ga_0}$ be any eigenfunction corresponding to $\al_1$. We can verify that $|u|\in \S\cap W_{\Ga_0}$ and $|u|$ is a minimiser of $J$ over $\S\cap W_{\Ga_0}$. Thus $|u|$ are also an  eigenfunction corresponding to $\al_1$. From the strong maximum principle either $|u|\equiv 0$ or $|u|>0$ in $\Om$. Since $u\neq 0$ we have $|u|>0$ in $\Om$, hence $\al_1$ is a principal eigenvalue of \eqref{Genp}.
\qed 

\begin{proposition}
	Let $p, \Om ,\al_1$ as in Proposition \ref{propo:appendix}. Let $u$ be an eigenfunction associated with $\al_1$. If $\Om\mbox{ and }\Ga_0$ are symmetric with respect to a hyperplane $H$, then $u$ is also symmetric with respect to the hyperplane $H$.
\end{proposition}
\textit{Proof.} Let $\sigma_H$ be the reflection with respect to the hyperplane $H$. Given  $\Om\mbox{ and }\Ga_0$ are symmetric with respect to a hyperplane $H$, so we have $x\in \Om\, (\mbox{or } \Ga_0) \text{ if and only if } \sigma_H(x)\in \Om \,(\mbox{or } \Ga_0).$ Consider an eigenfunction $u$ corresponding to $\al_1$ and define $v(x)\dq u(\sigma_H(x)) \text{ for } x\in\Om$. Then $v\in W_{\Ga_0}$ with $u=v$ on $H$ and $J(v)=J(u)=\al_1.$
By the simplicity of $\al_1$, we have $v=ku\text{ a.e. in } \Om$ for some $k>0$. Since $\|u\|_p=\|v\|_p$, we have $k=1.$ Therefore $v=u \text{ a.e. in } \Om$, i.e., $u(x)=u(x^*)\text{ a.e. } x\in \Om$. Hence $u$ is symmetric with respect to $H$.
\qed
\begin{remark}
 In particular if $\Om \mbox{ and } \Ga_0$ are radially symmetric, then the eigenfunction $u$ associated with $\al_1$ is radial in $\Om$. Thus the eigenfunctions associated with $\nu_1(\Om^\#)$ and $\tau_1(\Om_\#)$ are radial in $\Om^\#$ and $\Om_\#$ respectively. Notice that, the standard arguments using the schwartz symmetrization and P\'{o}lya-Szeg\"{o} inequality does not work, if the domain $\Om$ is not a ball. 
 
\end{remark}
\section{The existence of infinitely many eigenvalues:} Using standard variational methods, as in \cite{Garcia87} for Dirichlet eigenvalues using Krasnoselskii genus, we can obtain a set of critical values of $J$ on $\S\cap W_{\Ga_0}$. For a symmetric closed subset $\A\subset\S$, Krasnoselskii genus of $\A$ is defined as
\[\ga(\A)\dq \inf\{n\in\N : \exists \mbox{ an odd continuous map from }\A \mbox{ into }\R^n\setminus\{0 \} \}\]
with the convention that $\inf\{\emptyset\}=\infty.$ For $n\in\N$, let
\[\E_n\dq\{\A\subset\S: \A=\overline{\A},\, \A=-\A\mbox{ and }\ga(\A)\geq n \}, \]
\[\la_n\dq \inf\limits_{\A\in\E_n}\sup\limits_{u\in\A}J(u). \]
For each $n\in \N,$ let 
\[\F_n\dq\left\{h(S^{n-1}): h \mbox{ is an odd continuous map from }S^{n-1} \mbox{ to }\S\cap W_{\Ga_0} \right\}.\]  
An another set of critical values of $J$ can be obtain as follows: 
\[\la_n^*\dq \inf\limits_{\A\in\F_n}\sup\limits_{u\in\A}J(u).\] 
Since $\ga(S^{n-1})=n,$ and $\ga$ is invariant under odd homeomorphisms, we have $\F_n\subseteq\E_n$ and hence $\la_n\leq \la_n^*,\, \forall n\in\N.$ Then, with similar arguments as in Proposition 5.3 of \cite{Garcia87} and 
Theorem 5 of \cite{Drabek1}, we can prove the following:
\begin{proposition}\label{propo:appendix 1}
Let $p\in(1,\infty)\mbox{ and } \Om$ as above. Then for $n\in\N,\, \la_n,\, \la_n^*$ are eigenvalues of \eqref{Genp} such that $\la_n,\,\la_n^*\nearrow\infty,\mbox{ as } n\rightarrow\infty.$
\end{proposition}

\section{The variational characterisation of second eigenvalue:}
We have that $\al_i=\la_i=\la_i^*,\, i=1,2.$ This follows for $i=1$, since the set $\{u,-u\}$ belongs to both $\E_1$ and $\F_1$ for $u\in\S\cap W_{\Ga_0}$. For $i=2$, this follows from the facts that the set $\{au_2^++bu_2^- : |a|^p\|u_2^+\|^p+|b|^p\|u_2\|^p=1 \}$ belongs to both $\E_2$ and $\F_2$, where $u_2$ is the second eigenfunction, and $\al_1$ is isolated (see Section 2 of \cite{Anoop-Drabek} for a similar characterisation of second Dirichlet eigenvalue). Now we have the following proposition.
\begin{proposition} Let $p\in (1,\infty)\mbox{ and }\Om$ be as in Theorem \ref{thm 1}. Let $\nu_2,\, \tau_2$ be the second eigenvalues of \eqref{ND} and \eqref{DN} respectively. Then  $\nu_2$ and  $\tau_2$ have the following variational characterisation:
\[\nu_2=\inf\limits_{\A\in\F_2}\sup\limits_{u\in\A} J(u),\ \tau_2=\inf\limits_{\A\in\G_2}\sup\limits_{u\in\A} J(u), \] where
$\G_2\dq\left\{h(S^{1}): h \mbox{ is an odd continuous map from }S^{1} \mbox{ to }\S\cap W_{\Ga_1} \right\}.$
\end{proposition}
\begin{acknowledgements}
The authors would like to thank Prof. S. Kesavan and Dr. Vladimir Bobkov for their valuable suggestions during the discussions, which helped to improve the manuscript. The first author  would like to  thank the Department of Science \& Technology, India for the research grant DST/INSPIRE/04/2014/001865.
\end{acknowledgements}
\bibliographystyle{spmpsci}      

\end{document}